\documentclass[11pt,twoside]{article}
\usepackage{verbatim,epsfig}
\usepackage{amsfonts,eucal}
\input xy
\xyoption{all}

\begin{document}

\setlength{\topmargin}{0in}
\setlength{\headheight}{9pt}
\setlength{\textwidth}{5in}
\setlength{\textheight}{8in}
\setlength{\paperheight}{11in}
\setlength{\paperwidth}{8.5in}
\setlength{\parskip}{2ex}

%
\newcommand{\barefootnote}[1] {%
  \begingroup
    \renewcommand{\thefootnote}{}
    \footnotetext{#1}
    \renewcommand{\thefootnote}{\arabic{footnote}}
  \endgroup
}

%
\renewcommand{\title}[1] {%
  \begingroup
    \begin{center}
      \vspace{0.4in}
      \bf\huge
      \addtolength{\baselineskip}{5mm}
      #1
    \end{center}
  \endgroup
}

\newcommand{\url}[1] {%
  \barefootnote{%
    {\small e-print archive: }
    {\texttt http://xxx.lanl.gov/#1}
  }
}

\renewcommand{\author}[1] {%
  \begingroup
    \begin{center}
      \vspace{0.4in}
      \bf
      #1
      \vspace{0.2in}
    \end{center}
  \endgroup
}

%
\newcommand{\address}[1] {%
  \begingroup
    \begin{center}
      #1
    \end{center}
  \endgroup
}

%
\newcommand{\addressemail}[1] {%
  \begingroup
    \begin{center}
      \vskip-\baselineskip
      #1
    \end{center}
  \endgroup
}

%
\newcounter{secondpage}
\newcommand{\cutpage}{%
  \newpage
  \setcounter{secondpage}{\value{gfirstpage}}
  \addtocounter{secondpage}{1}
  \setcounter{page}{\value{secondpage}}
}

%
\newcounter{gfirstpage}
\newcounter{glastpage}
\newcommand{\copyrightnotice}[4]{
  \leftline{\copyright~ #1 International Press}
  \leftline{Adv. Theor. Math. Phys. {\bf #2} (#1) #3--#4}
  \setcounter{gfirstpage}{#3}
  \setcounter{glastpage}{#4}
}


\pagestyle{myheadings}
\thispagestyle{empty}

\newcommand{\Diff}{{\rm Diff_{+}}}
\newcommand{\Mor}{{\rm Mor}}
\newcommand{\C}{{\mathbb C}}
\newcommand{\Q}{{\mathbb Q}}
\newcommand{\R}{{\mathbb R}}
\newcommand{\Z}{{\mathbb Z}}
\newcommand{\HF}{{\bf HF}}
\newcommand{\tg}{{\tau}}
\newcommand{\kk}{{\bf k}}
\newcommand{\tr}{{\rm Trace \;}}
\newcommand{\Gras}{{\rm Grass}}
\newcommand{\M}{{\bf M}}
\newcommand{\Met}{{\rm Metrics}}
\newcommand{\Aut}{{\rm Aut}}
\newcommand{\Oh}{{\rm O}}
\newcommand{\SO}{{\rm SO}}
\newcommand{\Spin}{{\rm Spin}}
\newcommand{\s}{{\mathbb S}}
\newcommand{\asd}{{\bf D}}
\newcommand{\ASD}{{\mathcal D}}
\newcommand{\Sym}{{\rm Sym}}
\newcommand{\Hom}{{\rm Hom}}
\newcommand{\SU}{{\rm SU}}
\newcommand{\U}{{\rm U}}
\newcommand{\SP}{{\rm SP}}
\newcommand{\wall}{{\rm Wall}}
\newcommand{\Bun}{{\rm Bun}}
\newcommand{\pt}{{\rm pt}}
\newcommand{\threehalf}{{\textstyle{\frac{3}{2}}}}

\title{A rudimentary theory of topological 4D gravity}
\url{math.DG/0007018}		
\author{Jack Morava}		
\address{Department of Mathematicss\\ The Johns Hopkins University\\
Baltimore 21218 Md.}
\addressemail{jack@math.jhu.edu}	
\markboth{\it PRETTY GOOD GRAVITY}{\it Morava}

\begin{abstract} A theory of topological gravity is a homotopy-theoretic 
representation of the Segal-Tillmann topologification of a two-category 
with cobordisms as morphisms. This note describes some relatively accessible 
examples of such a thing, suggested by the wall-crossing formulas of 
Donaldson theory. \end{abstract}

\section {Gravity categories}

\noindent
A {\bf cobordism category} has manifolds as objects,
and cobordisms as morphisms. Such categories were introduced by
Milnor [22], but following Segal's definition of conformal field 
theory [29] and Atiyah's subsequent abstraction of the notion of 
topological quantum field theory [1] they have been studied very
widely. Recently, Tillmann [31] has shown the utility in this context
of certain closely related {\bf two}-categories (which generalize
the classical notion of category, by admitting morphism-objects
which are themselves categories). The following definition is based 
on her ideas. 

\cutpage

\noindent
{\bf Definition} A {\bf gravity} two-category has \medskip

\noindent
$\bullet$ (closed) {\bf manifolds} as objects, \medskip
  
\noindent
$\bullet$ {\bf cobordisms} as morphisms, and \medskip

\noindent
$\bullet$ {\bf isomorphisms} of these cobordisms, equal to the
identity on the boundary, as {\bf two}-morphisms. \bigskip

\noindent
There are many possible variations on this theme, and I will not
try for maximal generality. If the objects of the category have
dimension $d$ (so the cobordisms are $(d+1)$-dimensional) then I
will say that the gravity category is $(d+1)$-dimensional. 
I will assume that manifolds are smooth, compact and oriented, 
but not necessarily connected, and (following Segal) I understand 
the empty set to be a manifold of any dimension.
\bigskip

\noindent
{\bf 1.1} If $V$ and $V'$ are $d$-manifolds, a morphism 
\[
   W : V \to V'
\]
is (the germ of) an orientation-preserving diffeomorphism
\[
   (V_{op} \cup V') \times [0,1) \cong \nu(\partial W)
\]
of the manifold on the left with a collar neighborhood of the
boundary of the $(d+1)$-manifold $W$; the subscript $op$ signifies 
reversed orientation. The morphism category $Mor(V,V')$ has such
cobordisms as its objects; it is a topological category, in which
the space of morphisms between two cobordisms $W$ and $\tilde W$ 
consists of orientation- and boundary-identification-preserving 
diffeomorphisms $W \cong \tilde W$. Gluing along the boundary 
defines a continuous composition functor
\[
W, W' \mapsto W \circ W' : Mor(V,V') \times Mor(V',V'') \to
Mor(V,V'') \;,
\]
while disjoint union of objects gives this two-category a
monoidal structure, with the empty set as identity object.
\bigskip

\noindent
By replacing $Mor(V,V')$ with its set $\pi_0
Mor(V,V')$ of equivalence classes of objects, we obtain the
category employed by Atiyah to define a topological quantum field
theory; in other words, we can pass from a gravity two-category,
in which the morphism objects are enriched by a categorical
structure, to a classical category, in which the morphism objects
are simply sets. Tillmann's more perspicacious alternative is to
interpret the topological category $Mor(V,V')$ as a simplicial
topological space and to replace it with its geometric
realization $\Mor(V,V')$. This construction preserves Cartesian
products (as does $\pi_0$: indeed the set of equivalence classes
of objects in $Mor$ is the set of components of the space
$\Mor$), defining a {\bf topological} gravity category (i.e., a
category in which the morphism objects are topological spaces,
and the composition maps are continuous). A topological quantum
field theory in the sense of Atiyah [12 \S 1.7] is thus a (continous)
monoidal functor from a topological gravity category to the
(topological) category of modules over a {\bf discrete}
topological ring. \bigskip  

\noindent
However, we can consider monoidal functors to more general
categories: for example, the singular chains on the morphism
spaces of a gravity category define a monoidal category enriched
over chain complexes, whose representations are the
(co)homological field theories of physics. In the language of
homotopy theory, these are representations in a category of
modules over some Eilenberg-MacLane ring-spectrum. In general, I
will call any monoidal functor from a topological gravity category to
the category of dualizeable objects over a ring-spectrum, a {\bf 
theory of topological gravity}. One of the points of this
paper is that there is a rich supply of such things. \bigskip

\noindent
{\bf 1.2} This terminology needs some explanation. If $W$ is a
manifold with boundary, let $\Diff(W)$ be the topological group
of orientation-preserving diffeomorphisms of $W$ which restrict
to the identity in some neighborhood of $\partial W$. The components 
of $\Mor(V,V')$ are indexed by equivalence classes of cobordisms 
$W : V \to V'$, and the components themselves are the classifying spaces 
$B\Diff(W)$. Gluing [20] defines a continuous homomorphism
\[
\Diff(W) \times \Diff(W') \to \Diff(W \circ W') \;;
\]
thus the (components of the) composition map in the topological
gravity category are the maps these compositions induce on 
classifying spaces. \bigskip

\noindent
On the other hand, a fundamental tautology of Riemannian geometry
asserts that an isometry of a complete connected Riemannian
manifold which fixes a frame at some point is the identity: such a
map preserves the geodesics out of the framed point, and any
other point in the manifold can be reached by such a geodesic. It
follows that group of diffeomorphisms framing some basepoint will
act {\bf freely} on the (contractible) space of Riemannian
metrics on a compact connected manifold. The space $B\Diff(W)$ is
the homotopy quotient of the space of metrics [10, 11] on $W$ by the 
diffeomorphism group and we can think of morphisms in
the $(d+1)$-dimensional gravity category as cobordisms between
$d$-manifolds, together with a choice of equivalence class of
Riemannian metric on the cobordism. Riemannian geometry thus
provides the gravity category with a smooth structure. \bigskip

\noindent
A (projective) Hilbert-space representation of a
topological gravity category, along the lines considered by Segal
in his definition of a conformal field theory, is thus very
close to a quantum theory of gravity. When $d = 1$ we can see
this more explicitly: the Riemann moduli space is the quotient of
the space of conformal structures on a closed connected surface
by the group of its orientation-preserving diffeomorphisms, which
acts with finite isotropy when the genus exceeds one. This
defines a monoidal functor from the two-dimensional gravity
category to Segal's, which (away from closed surfaces of low
genus) is a rational homology isomorphism on morphism spaces.
Consequently, any conformal field theory in Segal's sense
defines a quantum theory of two-dimensional gravity. \bigskip   

\noindent
{\bf 1.3 Examples:} \medskip

\noindent
i) From this point of view, there is no {\it a priori} reason to 
limit ourselves to smooth manifolds. We could begin with a two-category 
of topological manifolds, and replace its morphism categories
by their classifying spaces, as before: there are plenty of 
non-smoothable four-manifolds! \medskip

\noindent
ii) In higher dimensions, the category of manifolds and
equivalence classes of $h$-cobordisms is a groupoid, with the
Whitehead group of an object as its automorphisms. In low
dimensions these categories are still quite mysterious. 

\noindent
iii) We can consider classes of manifolds with extra structure: 
for example, by requiring that the Stiefel-Whitney class $w_2$ vanish, we 
can define a gravity category of four-dimensional Spin-manifolds. 
[The set of Spin-structures on such a manifold is a principal homogeneous
space over its first mod two cohomology group, but is not naturally
isomorphic to that group.] \medskip

\noindent
iv) Similarly, the four-dimensional gravity category of $\Spin^{\C}$-manifolds 
is defined by cobordisms endowed with a complex line bundle with Chern 
class lifting $w_2$. \bigskip

\noindent
Ex.\ iii) can be regarded as the subcategory of Ex.\ iv) defined by objects
with trivial Chern class. It is natural to think of the morphism categories
in Ex.\ iii) as graded by elements of the middle homology lattice; for 
example, algebraic surfaces lie on the quadric $c_1^2 = 2 \chi + 3 \sigma$.
[Note that reversing orientation changes the signature, but not the
Euler characteristic.] \bigskip   

\noindent
When $d$ is {\it odd}, the morphisms of a $d+1$-dimensional gravity
category are naturally graded by Euler characteristic: the correction
term in the formula
\[
\chi(W \circ W') = \chi(W) + \chi(W') - \chi(W \cap W')
\]
is zero. When $d$ is one, the Euler characteristic counts the number of 
handles or loops in the usual quantum or genus expansion; it defines a zeroth 
Mumford class $\kappa_0$. If we exclude closed manifolds from our morphism 
spaces, and thus do not admit the empty set as a plausible object, this 
grading is bounded below. The signature defines a similar grading, when
$d=3$. \bigskip

\noindent
Many interesting decorations of gravity categories are possible: Lorentz 
cobordism [28, 33], defined by a nowhere-vanishing vector field oriented 
suitably at the boundary, is one example. Restricting the objects (e.g. to 
be unions of (standard, or homology) spheres, or contact manifolds 
[19]) is another alternative. Witten's original two-dimensional theory [34] 
admits singular (stable) algebraic curves as morphisms; this compactifies its 
morphism spaces, and Kontsevich has shown (as Witten conjectured) that the 
resulting theory has a well-behaved vacuum state. \bigskip 

\section {Pretty good theories of topological gravity}

\noindent
A Riemannian metric $g$ on an oriented closed connected 
two-manifold $\Sigma$ defines a Hodge operator $*_g$ on its harmonic 
forms. This operator squares to $-1$ on one-forms, and so defines a 
complex structure on the de Rham cohomology $H^1_{dR}(\Sigma)$. The 
space of isomorphism classes of complex structures on a real Euclidean 
space of dimension $2g$ is the quotient $\SO(2g)/\U(g)$, so we get a map
\[
\tau : B\Diff(\Sigma) \to (\Met)/(\Diff(\Sigma)) \to \SO/\U 
\]    
in the large genus limit. This can be constructed more generally 
by working with differential forms which vanish on the boundary. 
Orthogonal sum of vector spaces makes an $H$-space of the target of 
$\tau$, and it is not hard to see that if $\Sigma$ and
$\Sigma'$ are surfaces with geodesic boundaries, then gluing them $c$ times
along some sets of compatible boundary components defines a 
homotopy-commutative diagram 
$$\xymatrix{
{B\Diff(\Sigma) \times B\Diff(\Sigma')} \ar[r] \ar[d]^{\tau \times \tau} & 
{B\Diff(\Sigma \circ \Sigma')} \ar[d]^{\tau}  \\
{\SO/\U \times \SO/\U} \ar[r]^{\oplus} & {\SO/U} \;.}
$$ 
[The intersection form on the middle homology of $\Sigma \circ \Sigma'$ is
the direct sum of the intersection forms of $\Sigma$ and $\Sigma'$, together
with a {\bf split hyperbolic} intersection form of rank $c-1$, which has
a canonical complex structure [32 IV \S 4].] \bigskip

\noindent
This is perhaps the simplest example of a theory of two-dimensional 
topological gravity: it is a monoidal homotopy-functor to a topological 
category with one object and the $H$-space $\SO/\U$ of morphisms 
[25]. The functor is a version of the Jacobian, which refines the infinite 
symmetric product construction (which takes disjoint union to Cartesian 
product). The Siegel moduli space for abelian varieties has the rational 
cohomology of an integral symplectic group, and a version of Hirzebruch's
proportionality principle implies that the stable rational cohomology of 
this moduli space agrees with the cohomology of $\SO/\U$. \bigskip

\noindent
{\bf 2.1.1} In general, a topological quantum field theory $\HF$ (with values 
in some category of modules over a ringspectrum $\kk$) assigns to a suitable 
$d$-manifold $V$, a module-spectrum $\HF(V)$, such that \medskip

\noindent
i) the construction is exponential, in the sense that 
\[ 
\HF(V \sqcup V') \cong \HF(V) \wedge \HF(V') \;;
\] 
\noindent
ii) there is a pairing 
\[
\tr: \HF(V_{op}) \wedge \HF(V) \to \kk
\]
which is nondegenerate, in the sense that the induced map from
$\HF(V_{op})$ to the functional $\kk$-dual of $\HF(V)$ is an isomorphism;
\medskip

\noindent
iii) there is a natural transformation 
\[
\tg_W : B\Diff(W) \to \HF(\partial W)
\]
\noindent
subject to a {\bf monoidal axiom}: if $\partial W = V_{op} \sqcup V'$, etc., 
then the diagram 
$$\xymatrix{
{B\Diff(W) \times B\Diff(W')} \ar[r] \ar[d]^{\tg \times \tg} & 
{B\Diff(W \circ W')} \ar[d]^{\tg}  \\
{\HF(V_{op}) \wedge (\HF(V') \wedge \HF(V_{op}')) \wedge \HF(V'')} 
\ar[r] & {\HF(V_{op}) \wedge \HF(V'')} \;.}
$$ 
commutes up to homotopy. \bigskip

\noindent
The smash product of two such functors yields another. \bigskip

\noindent
{\bf 2.1.2} Objects in the two-dimensional gravity category are just 
collections of circles, which can be indexed by nonnegative integers. In 
this case, a theory is defined by a dualizable $\kk$-module spectrum $\M$, 
together with a system
\[
\tg^p_q \in (M^{\wedge (p+q)})^*(B\Diff (\Sigma)) = [B\Diff (\Sigma),
\M \wedge_{\kk} \dots \wedge_{\kk} \M]^*
\]
of characteristic classes for bundles of connected surfaces $\Sigma$ with 
$p$ incoming and $q$ outgoing boundary components, which behave compatibly 
under gluing. The example above is deceptive, for in that case $\M$ agrees 
with the group ring $\kk = \s[\SO/\U]$, so the multiple smash product 
simplifies. The topological category with one object, and Tillmann's 
group-completion
\[
\coprod_{g \geq 0} B\Diff(\Sigma_g) \to \Z \times B\Gamma^+_{\infty}
\]
as its space of morphisms, defines the universal example of a theory of
this type; the cohomology homomorphism defined by the induced map
\[
\Z \times B\Gamma^+_{\infty} \to \SO/\U
\]
factors through the classical map which kills the Mumford classes in
degree divisible by four. In more general cases related to quantum 
cohomology [20, 24], $\M$ will be a Frobenius object in the category of 
spectra, and the theory can be reformulated in terms of a family of natural 
transformations
\[
\otimes^{p+q}H^*(\M) \to H^*(B\Diff(\Sigma)) \;.
\] \medskip
 
\noindent
{\bf 2.2} The Hodge-theoretic construction described above has a close 
analogue for four-manifolds, which is also classical in a way: the 
wall-crossing formulas [17] of Donaldson theory are its descendants. As 
in dimension two, its construction is based on properties of the intersection 
form on middle cohomology: \bigskip 

\noindent
If $W$ is an compact connected oriented four-manifold with 
$\partial W$ a union of homology spheres then the intersection form
\[
   x,y \mapsto \langle x, y \rangle = (x \cup y)[W,\partial W]
\] 
on the integral lattice $B = H^2(W,\partial W,\Z)$ is unimodular. In 
dimension four, Wu's formula implies that
\[
   q(x) = \langle x,x \rangle \equiv \langle x, w_2 \rangle
\]
modulo two, so the form $q$ is even if the manifold admits a Spin-structure
[16 \S 5.7.6]. On a $\Spin^{\C}$-manifold the intersection form is even or odd 
depending on the parity of the Chern class of its associated complex line 
bundle. \bigskip

\noindent
By a fundamental theorem of Freedman [13] any unimodular quadratic
form can arise as the intersection form of a closed topological four-manifold;
but by similarly fundamental results of Donaldson [6, 9] the intersection form 
of a closed {\bf smooth} four-manifold is either indefinite, or diagonalizable 
over the integers.  As in two dimensions, the action of a diffeomorphism
on homology defines a monodromy representation
\[
   \Diff(W) \to \Aut_+(B,q) = \SO(B) 
\]
which factors through $\pi_0(\Diff(W))$; it is convenient to think of
its kernel [18] as an analogue, for four-manifolds, of the Torelli group of
surface theory. \bigskip

\noindent
{\bf 2.3} Let $b = b_+ + b_-$ be the rank, and $\sigma = b_+ - b_-$ the 
signature, of the inner product space defined by $q$ on $B \otimes \R$. For 
our purposes the {\bf indefinite} lattices are the most interesting: these are 
classified by their rank, signature, and type (even if $q(x) \equiv 0$ mod two,
otherwise odd). In the indefinite case, the manifold $\Gras^-(B)$ of maximal 
negative-definite subspaces of $B \otimes \R$ is a noncompact (contractible) 
symmetric space defined by a cell of dimension $b_+b_-$ in the usual 
Grassmannian of $b_-$-planes in $b$-space. The orthogonal group of the lattice 
acts on this cell with finite isotropy, so the canonical homotopy-to-geometric 
quotient map
\[
   B\SO(B) \to \Gras^-(B)/\SO(B)
\]
is a rational homology isomorphism. If $B$ and $B'$ are indefinite lattices,
then the construction which sends a pair of negative definite subspaces in 
the real span of each, to their orthogonal sum in the real span of the direct 
sum lattice, defines a map
\[
\Gras^-(B) \times \Gras^-(B') \to \Gras^-(B \oplus B')
\]
which is equivariant with respect to the Whitney sum homomorphism
\[
\SO(B) \times \SO(B') \to \SO(B \oplus B')
\]
The Grothendieck group of the category of even indefinite unimodular lattices 
is free abelian on two generators, corresponding to the hyperbolic plane and 
the $E_8$ lattice [30 V \S 2]. The `Hasse-Minkowski' spectrum ${\bf K}_{\rm EIU}$ 
defined by the algebraic $K$-theory of the category of such 
lattices is the group completion of the monoid constructed from the disjoint 
union of the classifying spaces of their orthogonal groups; the tensor product 
of two such lattices defines another, making this a commutative ring-spectrum. 
\bigskip

\noindent 
{\bf 2.4} A Riemannian metric $g$ on $W$ defines a Hodge operator $*_g$ on 
harmonic forms, but now this operator squares to $+1$ on the middle cohomology.
The function which assigns to $g$, the $*_g = -1$-eigenspace of harmonic 
two-forms vanishing on $\partial W$, maps the space of Riemannian metrics 
to the negative-definite Grassmannian $\Gras^-(B)$ equivariantly
with respect to the action of $\Diff(W)$. \bigskip

\noindent
If $W$ and $W'$ are four-manifolds bounded (as above) by homology spheres,
and if $W \circ W'$ results from gluing these manifolds along a collection of
compatible boundary components, then the quadratic module of $W \circ W'$ is
canonically isomorphic to $B \oplus B'$; hence the cohomology representation
of the diffeomorphism group defines a monoidal functor from the gravity 
category of Spin four-manifolds bounded by standard spheres, to the topological
category with one object, and the Hasse-Minkowski spectrum as morphisms. There
is a similar functor defined on the category with homology spheres as objects,
but the resulting lattice is no longer necessarily indefinite [9 \S 1.2.3]. 
\bigskip

\noindent
The higher algebraic $K$-theory of such lattices has apparently not 
received much attention. It is remarkable that the relatively naive
constructions sketched above already define pretty good theories of 
topological gravity. The $\eta$-invariant of Atiyah-Patodi-Singer [3]
is much more sophisticated; to find an interpretation in these terms
for it, analogous to the way Floer homology globalizes the Casson 
invariant, would be extremely interesting. \bigskip

\section {Toward a parametrized Donaldson theory}

\noindent
A good theory of gravity shouldn't exist in a vacuum: it deserves to be 
coupled to some nontrivial matter. Donaldson [8] and Moore and Witten 
[23] have suggested the study of equivariant supersymmetric Yang-Mills theory 
parameterized by classifying spaces of diffeomorphism groups. A
fragment of such a theory is sketched below. \bigskip

\noindent 
{\bf 3.1} Suppose for simplicity that $W$ is closed. The graded space 
$\Bun_*(W)$ of gauge equivalence classes of connections on $\SU(2)$-bundles 
over $W$ has components indexed by the second Chern class of the bundle. 
Let $\asd_*$ be the subspace of $\Met \times \Bun_*(W)$ consisting of pairs 
$(g,A)$, where $A$ is a connection on an $\SU(2)$-bundle over $W$ with 
curvature two-form
\[
   *_g(F_A) = - F_A 
\]
antiselfdual with respect to the metric $g$. The standard transversality 
arguments of Donaldson theory [9 \S 4.3] imply that this space is a manifold, 
with fiber of dimension $8c_2 - \threehalf(\sigma + \chi)$ above the metric 
$g$; at least, provided this metric admits no {\bf reducible} antiselfdual 
connections. These reducible connections define an interesting kind of 
distinguished boundary for the space of antiselfdual connections. \bigskip

\noindent
{\bf 3.2} Reducible connections on $W$ are parametrized by the wall arrangement
\[
\wall(B) = \{ H \in \Gras^-(B) \;|\; H \cap B \neq \{ 0 \} \; \}
\]
of the lattice $B$: it is the set of maximal negative-definite subspaces
of $B \otimes \R$ containing a lattice point. This is a union of smooth
submanifolds of codimension $b_-$, filtered by the increasing 
family $\wall_d(B)$ of subspaces consisting of maximal negative-definite $H$
containing a lattice point $x$ with $0 > q(x) \geq - d$ (which is a locally
finite union of manifolds [14]). The orthogonal group of $B$ acts naturally
on these arrangements, as well as on the quotient spaces 
\[
  \wall^S_d(B) = \Gras^-(B)/\wall_d(B)
\]
(which are roughly the $S$-duals of the wall arrangements). If $B$ and $B'$ 
are two indefinite lattices, then the orthogonal direct sum map defines a 
commutative diagram 
$$\xymatrix{
{\Gras^-(B) \times \Gras^-(B')} \ar[r] \ar[d] & 
{\Gras^-(B \oplus B')} \ar[d] \\ {\wall^S_d(B) \wedge \wall^S_{d'}(B')} 
\ar[r] & {\wall^S_{d+d'}(B \oplus B')}}
$$ 
which is equivariant, with respect to the Whitney sum on orthogonal groups. 
The equivariant cohomology $H^*_{\SO}(\wall^S_*)$ defines yet another variant 
of a topological gravity theory, but there seems to be little known about 
such essentially arithmetic invariants. \bigskip

\noindent
{\bf 3.3} If $g$ is in the complement of the preimage $\Met^0_d$ of 
$\wall_d$ in the space $\Met$ of metrics on $W$, then no $\SU(2)$-bundle
with Chern class less than $-d$ admits a connection with $*_g$-antiselfdual 
curvature. Thus if $\asd^0_d$ denotes the space of pairs $(g,A)$ such that
$A$ is gauge equivalent to a connection induced from a line bundle with
curvature antiselfdual with respect to $g$, then
\[
  (\asd_d,\asd^0_d) \to (\Met,\Met^0_d) \times \Bun_d(W)
\]
is a kind of $\Diff(W)$-equivariant cycle, of relative finite dimension above
the space of metrics. It cannot be expected to be proper, but Donaldson 
theory has developed sophisticated methods to deal with such issues [7]: 
let $\SP^{\infty}_d(W_+)$ be the space of finitely supported functions $f$ 
from $W$ to the integers, such that 
\[
   \sum_{x \in W} f(x) = d \;,
\] 
and let 
\[
\overline {\asd}_d = \coprod_{0 \leq i \leq d} \asd_i \times 
\SP^{\infty}_{d-i}(W_+)
\]
be the analogue of the Uhlenbeck-Donaldson compactification of $\asd_d$
in the stratified space
\[
\Met \times (\coprod_{0 \leq i \leq d} \Bun_i(W) \times 
\SP^{\infty}_{d-i}(W_+)) = \Met \times \overline{\Bun}_d(W) \;. 
\]
Completing the subspace $\asd^0_d$ of reducible connections analogously 
defines a candidate
\[
(\overline{\asd}_d,\overline{\asd}^0_d) \to (\Met,\Met^0_d) \times 
\overline{\Bun}_d(W)
\]
for a $\Diff(W)$-equivariant Donaldson cycle. \bigskip

\noindent
To extract homological information from this construction, note that 
a class $z$ of dimension $*$ in the rational homology of $B\Diff(W)$ 
maps to a sum, with rational coefficients, of homology classes defined 
by maps 
\[
   Z \to \Met \times_{\Diff} \pt
\]
of {\bf smooth} manifolds $Z$. The fiber product of such a map with the projection
\[
\overline{\asd}_d \to \Met \times_{\Diff}\overline{\Bun}_d(W) \to \Met 
\times_{\Diff} \pt
\]
defines a class of dimension $* + 8d - \threehalf (\sigma + \chi)$ in the 
rational homology of 
\[
  (\Met,\Met^0_d) \times_{\Diff} \overline{\Bun}_d(W) \;;
\]
note that this admits a canonical map to the space
\[
\wall^S_d \wedge_{\SO(B)} \SP^{\infty}_d(W_+) \;,
\]
which depends only on the lattice $B$. \bigskip

\noindent
{\bf 3.4} The homotopy-to-geometric quotient map for the space of connections 
is a rational homology equivalence of $\Bun_*(W)$ with the space of based 
smooth maps from $W_+$ to $B\SU(2)$ [9 \S 5.1.15], and the Pontrjagin class 
defines a rational homology isomorphism of the space of maps with the 
Eilenberg-MacLane space $H(\Z,4)$. By the Dold-Thom theorem,
\[
\pi_i{\rm Maps}(W_+,H(\Z,4)) \cong H^{4-i}(W,\Z) \cong H_i(W,\Z) \cong 
\pi_i(\SP^{\infty}(W_+)) 
\]
so as far as rational (co)homology is concerned, we can replace the space 
$\overline{\Bun}_*(W)$ with the free topological abelian group on $W$. 
[This identification uses Poincar\'e duality, and hence requires a choice of 
orientation: the space of bundles is a contravariant functor, but the infinite 
symmetric product is covariant.] Combined with the constructions outlined 
above, this defines a generalized Donaldson invariant as a homomorphism
\[
\ASD_d : H_*(B\Diff,\Q)) \to H_{* + 8d - \threehalf (\sigma + \chi)}(\wall^S_d 
\wedge_{\SO} \SP^{\infty}_d,\Q)
\]
with values in a group which depends only on the cohomology lattice $B$; indeed
the rational homology of $\SP^{\infty}(W_+)$ is the symmetric algebra on
the homology of $W$, and the automorphic cohomology 
\[
  H^*_{\SO(B)}(\SP^{\infty}(W_+),\Q) = H^*(\SO(B),{\rm Sym}(H^*(W)))
\]
contains the classical ring of automorphic forms for the orthogonal
group [5] as the invariant elements of the symmetric algebra on $B$.
\bigskip

\noindent
This invariant generalizes the usual one, in the sense that $\ASD_d$
on a generator of the zero-dimensional homology of $B\Diff$  is the classical
invariant. [The usual convention is to interpret the antiselfdual 
cycle as a function on the {\bf co}homology of $W$, by taking its 
Kronecker product with $\exp(x), x \in H^*(W)$.] A four-manifold is 
said to be of {\bf simple} type, if the behavior of its classical 
invariant as a function of charge is not too complicated: in the present 
formalism, the condition is that
\[
\ASD_{d+1}(1) \mapsto w_0w_4^2 \ASD_d(1)
\]
under the homomorphism induced by the restriction map from $\wall^S_{d+1}$ to 
$\wall^S_d$ (where $w_0$ and $w_4$ generate the homology in degrees zero and 
four of $W$). This suggests 
\[
\tilde \ASD_d = (w_0w_4^2)^{-d} \ASD_d \in {\rm Hom}^{-\threehalf (\sigma + \chi)}
(H_*(B\Diff),H_*(\wall^S_d \wedge_{\SO} \SP^{\infty}_0))
\]
as the natural normalization for the generalized invariant.
\bigskip
 
\section{On the inadequacy of the foregoing}

\noindent
The preceding sketch defines at best a {\bf piece} of a topological 
gravity functor. It is defined only for manifolds without boundary, but it
behaves correctly under disjoint union: if $W_0$ and $W_1$ are two closed 
four-manifolds, then
\[
\sum_{d = d_0 + d_1} \ASD_{d_0}(W_0) \otimes \ASD_{d_1}(W_1) \mapsto
\ASD_d(W_0 \cup W_1)
\]
under the maps of \S 3.2; this is nothing but a definition of
the generalized invariant for non-connected manifolds. \bigskip

\noindent
In fact there is reason to think that these constructions may have
wider validity. Some years ago, Atiyah [2] proposed a unification 
of the invariants of Donaldson and Floer, based on a theory of 
semi-infinite cycles in the polarized manifold of connections on 
a three-manifold. A theory of such cycles which behaves naturally 
under variation of the metric on a bounding four-manifold would yield 
a topological gravity theory for four-manifolds, taking values in 
generalized automorphic forms with coefficients in Floer homology. \bigskip 

\noindent
Many results which follow from Atiyah's program are known now to be true; 
but (mostly because of difficulty with compactifications), work on these 
questions has advanced without using his cycle calculus. I am told, however, 
that recently there has been progress along the lines he suggested [26], 
though in Seiberg-Witten rather than Floer-Donaldson theory. Meanwhile,
Bauer [4] and Furuta [15] have studied generalized Seiberg-Witten 
invariants from a homotopy-theoretic point of view, and Bauer has shown that
his invariant behaves nicely under connected sum. The hope that these
new developments can be extended to the context proposed in this paper
has encouraged me to write this incomplete and probably naive account. 
\bigskip 

\noindent
{\bf Acknowledgements} This research was supported by the NSF. It is a 
pleasure to thank Stefan Bauer, Paul Feehan, Kenji Fukaya, Mikio Furuta, 
Peter Oszvath, Andrei Tjurin, and Richard Wentworth for helpful conversations 
about the material in this paper. It is, however, very speculative, and they 
deserve no blame for my excesses. \bigskip 

\bibliographystyle{amsplain}

\end{document}